\documentclass[11pt,leqno]{article}

\usepackage[english]{babel}

\usepackage{amsmath}

\usepackage{amssymb}

\usepackage{xypic}

\usepackage{amsfonts}

\usepackage{mathabx}

\numberwithin{equation}{section}

\newcommand{\A}{\mathcal{A}}

\newcommand{\C}{\mathcal{C}}

\newcommand{\D}{\mathcal{D}}

\renewcommand{\S}{\mathcal{S}}

\newcommand{\W}{\mathcal{W}}

\renewcommand{\mod}{\mathrm{Mod}}

\newcommand{\op}{\mathrm{Op}}

\newcommand{\rc}{\mathbb{R}\textrm{-}\mathrm{c}}

\newcommand{\CC}{\mathbb{C}}

\newcommand{\R}{\mathbb{R}}

\newcommand{\N}{\mathbb{N}}

\newcommand{\II}{\mathcal{I}}

\newcommand{\iso}{\stackrel{\sim}{\to}}

\newcommand{\CWM}{\C^{{\infty ,\mathrm{w}}}_M}

\newcommand{\wtens}{\overset{\mathrm{w}}{\otimes}}

\newcommand{\rh}{\mathit{R}\mathcal{H}\mathit{om}}

\newcommand{\ho}{\mathcal{H}\mathit{om}}

\newcommand{\Rh}{\mathrm{RHom}}

\newcommand{\id}{\mathrm{id}}

\renewcommand{\dim}{\textbf{Proof.}}

\newcommand{\qed}{\nopagebreak \phantom{} \hfill $\Box$ \\}

\newcommand{\imin}[1]{#1^{-1}}

\newcommand{\lind}[1]{\underset{#1}{\underrightarrow{\lim}}}

\newtheorem{teo}{Theorem}[section]%[chapter]

\newtheorem{df}[teo]{Definition}%[chapter]

\newtheorem{cor}[teo]{Corollary}%[chapter]

\newtheorem{oss}[teo]{Remark}%[chapter]

\newtheorem{lem}[teo]{Lemma}%[chapter]

\author{Luca Prelli}

\title{De Rham theorem for Whitney functions}

\date{}

\begin{document}

\maketitle

\begin{abstract}
Let $M$ be a real analytic manifold, $F$ a bounded complex of constructible sheaves. We show
that the Whitney-de Rham complex associated to $F$ is quasi-isomorphic to $F$.
\end{abstract}

\section{Introduction}

Let $M$ be a real analytic manifold, by Poincar\'e Lemma, it is well-known that
the de Rham complex with $\C^\infty$-coefficients over $M$ is isomorphic to $\CC_M$. In this paper we
show that $\CC_M$ is isomorphic to the de Rham complex with Whitney coefficients over $M$ on the associated subanalytic site.
Then, using sheaf theoretical arguments we can easily prove that, given a bounded complex of constructible sheaves $F$, the Whitney-de Rham complex associated to $F$ is quasi-isomorphic to $F$.
As a corollary we obtain a theorem of \cite{BP08}. (Another proof was given in \cite{Ch15} using deep results on $\D$-modules.) We also obtain a de Rham theorem for Schwartz functions on open subanalytic subsets. \\

\noindent Our proof uses these three known facts:
\begin{itemize}
\item There exists a site, the so called subanalytic site, where Whitney functions form a sheaf.
\item Locally on the site, the sections of this sheaf are nothing but $\C^\infty$-functions with bounded derivatives.
\item The homotopy axiom holds for de Rham cohomology with bounded derivatives.
\end{itemize}
Let us detail our arguments. \\

\noindent It is well known that Whitney $\C^\infty$-functions on the real analytic manifold $M$ do not satisfy the gluing conditions on open coverings. Hence, they do not form a sheaf. However, following \cite{KS01}, one can overcome this problem by associating to $M$ a Grothendieck topology by choosing as open subsets subanalytic open subsets and as coverings locally finite ones. This is the subanalytic site $M_{sa}$, and here Whitney $\C^\infty$-functions form a (subanalytic) sheaf denoted $\CWM$. \\

\noindent On $M_{sa}$, one can also define the (subanalytic) sheaf of $\C^\infty$-functions with bounded derivatives on bounded open subanalytic subsets, and prove that the associated de Rham complex is quasi-isomorphic to the constant sheaf on $M$. \\

\noindent $\C^\infty$-functions on $U$ with bounded derivatives coincide with Whitney $\C^\infty$-functions on $\overline{U}$ when $U$ is 1-regular (see \cite{Wh34}). Thanks to a decomposition result of \cite{Pa15} (which uses techniques developed in \cite{Pa94}) one can prove that locally on $M_{sa}$ the sheaf of bounded $\C^\infty$-functions coincide with $\CWM$. \\

\noindent Moreover, applying some (subanalytic) sheaf theoretical techniques one can link the de Rham complex with coefficients in $\CWM$ with the one with Whitney coefficients on any closed subanalytic subset $Z$ of $M$. This is because, basically, $\Rh(D'\CC_Z,\CWM)$ is quasi-isomorphic to the space of Whitney $\C^\infty$-functions on $Z$. Here $\CC_Z$ denotes the constant sheaf on $Z$ and $D'\CC_Z$ its dual $\rh(\CC_Z,\CC_M)$. This is done thanks to the Whitney tensor product introduced by Kashiwara-Schapira in \cite{KS96}. When $Z=\overline{U}$, with $U$ open 1-regular relatively compact and such that $D'\CC_U \simeq \CC_Z$, we obtain that $\Rh(D'\CC_Z,\CWM)$ is the space of $\C^\infty$-functions on $U$ with bounded derivatives. \\

\noindent Thanks to these considerations, one can reduce the isomorphism between the cohomology of $Z$ and the  Whitney de Rham complex to the quasi-isomorphism between $\CC_M$ and the de Rham complex with coefficients in $\CWM$, which is locally the homotopy axiom for de Rham cohomology with bounded derivatives and can be obtained as in \cite{BT82}. \\

\noindent The statement still holds if we replace $\CC_Z$ with a bounded complex of constructible sheaves $F$ and  Whitney $\C^\infty$-functions with the Whitney tensor product $F \wtens \C^\infty_M$. This allows to consider different kinds of de Rham theorems, as the one for Schwartz functions when $F=\CC_U$, $U$ open subanalytic.

\section{Subanalytic sheaves}

The following results on subanalytic sheaves are extracted from \cite{KS01} (see also \cite{Pr1}). We refer to \cite{KS90} for classical sheaf theory. \\

Let $X$ be a real analytic manifold. Denote
by $\op(X_{sa})$ the category of open subanalytic subsets of $X$. One
endows $\op(X_{sa})$ with the following topology: $S \subset
\op(X_{sa})$ is a covering of $U \in \op(X_{sa})$ if for any
compact $K$ of $X$ there exists a finite subset $S_0\subset S$
such that $K \cap \bigcup_{V \in S_0}V=K \cap U$. We will call
$X_{sa}$ the subanalytic site.\\

Let $\mod(\CC_{X_{sa}})$ (resp. $D^b(\CC_{X_{sa}})$) denote the category of sheaves (resp. bounded complexes of sheaves) of $\CC$-vector spaces on $X_{sa}$
and let $\mod_{\rc}(\CC_X)$  (resp. $D^b_{\rc}(\CC_X)$) be the abelian category of
$\R$-constructible sheaves (resp. bounded complexes of sheaves with $\R$-constructible cohomology) on $X$.\\

We denote by $\rho: X \to X_{sa}$ the natural morphism of sites.
We have functors
\begin{eqnarray*}
\rho_*:\mod(k_X) & \to & \mod(k_{X_{sa}}) \\
\imin \rho:\mod(k_{X_{sa}}) & \to & \mod(k_X) \\
\rho_!:\mod(k_X) & \to & \mod(k_{X_{sa}})
\end{eqnarray*}

%\begin{equation*}
%\xymatrix{\mod_{\rc}(\CC_X) \subset \mod(k_X)   \ar@ <2pt>
%[r]^{\hspace{1cm}\rho_*} &
%  \mod(k_{X_{sa}}) \ar@ <2pt> [l]^{\hspace{1cm}\imin \rho}. }
%\end{equation*}
The functors $\imin \rho$ and $\rho_*$ are the functors of inverse
image and direct image associated to $\rho$. The functor $\rho_!$ is left adjoint to $\imin\rho$.
It is fully faithful and exact. If $F \in \mod(\CC_{X_{sa}})$, $\rho_!F$ is the sheaf associated to the presheaf $U \mapsto \Gamma(\overline{U};F)$.
The functor $\rho_*$ is fully faithful and exact on
$\mod_{\rc}(\CC_X)$ and we identify $\mod_{\rc}(\CC_X)$ with its
image in $\mod(\CC_{X_{sa}})$ by $\rho_*$.

\section{Whitney tensor product}

Let $M$ be a real analytic manifold. We denote by $\C^\infty_M$ (resp. $\A_M$) the sheaf of $\C^\infty$-functions (resp. analytic functions) on $M$ and by $\D_M$ the sheaf of differential operators on $M$ with analytic coefficients. We denote by $\Omega^p_M$ the sheaf of $p$-differential forms with coefficients in $\A_M$ and by $\Omega^\bullet_M$ the complex
$$
0 \to \Omega^0_M \to \cdots \to \Omega_M^{{\rm dim}M} \to 0
$$
We denote by $\mod(\D_M)$ (resp. $D^b(\D_M)$) the category (resp. bounded derived category) of sheaves of $\D_M$-modules. References are made to \cite{KS96} for a complete exposition on formal cohomology.

\begin{df} Let $Z$ be a closed subset of $M$. We denote by
$\II^\infty_{M,Z}$ the sheaf of $\C^\infty$-functions on $M$
vanishing up to infinite order on $Z$.
\end{df}

\begin{df} A Whitney function on a closed subset $Z$ of $M$ is an indexed
family $F=(F^k)_{k\in \N^n}$ consisting of continuous functions on
$Z$ such that $\forall m \in \N$, $\forall k \in \N^n$, $|k| \leq
m$, $\forall x \in Z$, $\forall \varepsilon >0$ there exists a
neighborhood $U$ of $x$ such that $\forall y,z \in U \cap Z$

$$\left|F^k(z)-\sum_{|j+k|\leq m}{(z-y)^j \over
j!}F^{j+k}(y)\right| \leq \varepsilon d(y,z)^{m-|k|}.$$ We denote
by $W^\infty_{M,Z}$ the space of  Whitney $\C^\infty$-functions on
$Z$. We denote by $\W^\infty_{M,Z}$ the sheaf $U \mapsto
W_{U,U\cap Z}^\infty$.
\end{df}

In \cite{KS96} the authors defined the functor
$$\cdot \wtens \C^\infty_M:\mod_{\rc}(\CC_M) \to \mod(\D_M)$$
in the following way: let $U$ be a subanalytic open subset of $M$
and $Z=M \setminus U$. Then $\CC_U \wtens
\C^\infty_M=\II^\infty_{M,Z}$, and $\CC_Z \wtens
\C^\infty_M=\W^\infty_{M,Z}$. This functor is exact and extends as
a functor in the derived category, from $D^b_{\rc}(\CC_M)$ to
$D^b(\D_M)$. Moreover the sheaf $F \wtens \C^\infty_M$ is soft for
any $\R$-constructible sheaf $F$.\\

Let $F \in D^b_{\rc}(\CC_M)$, one sets
$$
F \wtens \C^{\infty,\bullet}_M := F \wtens \C^\infty_M \underset{\A_M}{\otimes} \Omega_M^\bullet.
$$

\section{The sheaf of Whitney functions}

Here we recall definition and some properties of the subanalytic sheaf of Whitney functions. References are made to \cite{KS01,Pr1}.

\begin{df} One denotes by $\CWM$ the sheaf on $M_{sa}$ of Whitney
$\C^\infty$-functions defined as follows:
$$U \mapsto \Gamma(M;H^0D'(\CC_U) \wtens \C^\infty_M),$$
where $D'(\cdot)=\rh(\cdot,\CC_M)$.
\end{df}

Remark that $\Gamma(U,\CWM)$ is a $\Gamma(\overline{U};\D_M)$-module for each $U \in \op(M_{sa})$, this implies that $\CWM$ has a structure of $\rho_!\D_M$-module.
We have the following result.
For each $F \in D^b_{\rc}(\CC_M)$ one has the isomorphism
$$
\imin \rho\rh(D'F,\CWM)  \simeq  F\wtens\C^\infty_M.
$$

Consider the complex $\Omega^\bullet_M$ of differential forms with analytic coefficients and set
$$
\C^{\infty,{\rm w},\bullet}_M := \CWM \underset{\rho_!\A_M}{\otimes} \rho_!\Omega_M^\bullet.
$$

One has the quasi-isomorphisms
$$
\imin\rho\rh(D'F,\C^{\infty,{\rm w},\bullet}_M) \simeq \imin\rho\rh(D'F,\CWM) \underset{\A_M}{\otimes} \Omega_M^\bullet \simeq F \wtens \C^{\infty,\bullet}_M
$$
where the first quasi-isomorphism follows from \cite[Proposition 2.2.3]{Pr1} and the fact that $\imin\rho \circ \rho_! \simeq \id$. \\

Let us consider a locally cohomologically trivial (l.c.t.) subanalytic open subset, i.e. $U \in \op({M_{sa}})$ satisfying $D'\CC_U \simeq \CC_{\overline{U}}$ and $D'\CC_{\overline{U}} \simeq \CC_U$. Then
$$
R\Gamma(U;\CWM) \simeq  R\Gamma(M;\CC_{\overline{U}}\wtens\C^\infty_M)
$$
is concentrated in degree zero. It means that $\CWM$ is the sheaf associated to the presheaf $U \mapsto \{\mbox{Whitney functions on $\overline{U}$}\}$. Indeed l.c.t. open subanalytic subsets form a basis for the topology of $M_{sa}$ (i.e. every locally compact subanalytic open subset has a finite cover consisting of l.c.t. subanalytic open subsets).

\section{Proof of the Theorem}

\begin{lem}\label{deRham} There is the following isomorphism in $D^b(\CC_{M_{sa}})$
$$
\CC_M \simeq \C^{\infty,{\rm w},\bullet}_M.
$$
\end{lem}
\dim\ \ It follows from \cite[Theorem 0.3]{Pa15} that every relatively compact subanalytic open subset has a finite cover consisting of contractible 1-regular open subanalytic subsets. A result of \cite{Wh34} asserts that on 1-regular open subsets Whitney functions are bounded $\C^{\infty}$-functions with bounded derivatives. Hence we are reduced to prove the isomorphism in $D^b(\CC_{M_{sa}})$
$$
R\Gamma(U;\CC_M) \simeq \Gamma(U;\C^{\infty,{\rm w},\bullet}_M).
$$
Since $U$ is contractible the left hand side is concentrated in degree 0 and $\Gamma(U;\CC_M) \simeq \CC$. We are reduced to prove that the right hand side is concentrated in degree zero and $H^0(U;\C^{\infty,{\rm w},\bullet}_M) \simeq \CC$. This is nothing but the homotopy axiom for de Rham cohomology with bounded derivatives, which can be obtained in the same way as the classical one, see \cite[Corollary 4.1.2]{BT82} (one just checks that the homotopy operator $K$ such that $f^*-g^*=K \circ d + d \circ K$ preserves bounded derivatives).
\qed

\begin{teo}\label{WDR} Let $F \in D^b_{\rc}(\CC_M)$. There is the following isomorphism in $D^b(\CC_M)$
$$
F \simeq F \wtens \C^{\infty,\bullet}_M.
$$
\end{teo}
\dim\ \ We have the chain of isomorphism in $D^b(\CC_M)$
\begin{eqnarray*}
F \wtens \C^{\infty,\bullet}_M & \simeq & \imin\rho\rh(D'F,\C^{\infty,{\rm w},\bullet}_M) \\
& \simeq & \rh(D'F,\CC_M) \\
& \simeq & F,
\end{eqnarray*}
where the second isomorphism follows from Lemma \ref{deRham} and the last one from the fact that $D'D'F \simeq F$ if $F \in D^b_{\rc}(\CC_M)$.
\qed

\begin{cor} \label{BP} Let $Z$ be a closed subanalytic subset of $M$. There is the following isomorphism in $D^b(\CC_M)$
$$
\CC_Z \simeq \mathcal{W}_{M,Z}^{\infty,\bullet}
$$
where $\mathcal{W}_{M,Z}^{\infty,\bullet}$ denotes the de Rham complex with coefficients in $\mathcal{W}_{M,Z}^{\infty,\bullet}$,
i.e. the Whitney-de Rham complex is isomorphic to $\CC_Z$.
\end{cor}
\dim\ \ Set $F=\CC_Z$ in Theorem \ref{WDR}. \qed

\begin{cor} \label{BPglobal} Let $Z$ be a closed subanalytic subset of $M$. There is the following isomorphism in $D^b(\CC)$
$$
R\Gamma(Z;\CC_Z) \simeq R\Gamma(M;\mathcal{W}_{M,Z}^{\infty,\bullet}),
$$
hence the cohomology $H^\bullet(Z)$ of $Z$ is isomorphic to the cohomology the de Rham complex with Whitney coefficients on $Z$.
\end{cor}
\dim\ \ We have the following chain of isomorphisms in $D^b(\CC)$
$$
R\Gamma(Z;\CC_Z) \simeq R\Gamma(M;\CC_Z) \simeq R\Gamma(M;\mathcal{W}_{M,Z}^{\infty,\bullet}),
$$
where the first isomorphism follows from the definition of direct image and the second one follows applying $R\Gamma(M;\cdot)$ to Corollary \ref{BP}.\qed

\begin{cor} \label{SDR} Let $U$ be an open subanalytic subset of $M$. There is the following isomorphism in $D^b(\CC)$
$$
R\Gamma_c(U;\CC_U) \simeq R\Gamma_c(M;\II^{\infty,\bullet}_{M,M \setminus U}),
$$
where $\II^{\infty,\bullet}_{M,M \setminus U}$ denotes the de Rham complex with coefficients in $\II^{\infty,\bullet}_{M,M \setminus U}$,
 hence the cohomology with compact support $H_c^\bullet(U)$ of $U$ is isomorphic to the cohomology with compact support of the de Rham complex with $\C^\infty$-coefficients on $M$ vanishing with all their derivatives outside $U$.
\end{cor}
\dim\ \ Set $F=\CC_U$ in Theorem \ref{WDR}. Then we have the following chain of isomorphisms in $D^b(\CC)$
$$
R\Gamma_c(U;\CC_U) \simeq R\Gamma_c(M;\CC_U) \simeq R\Gamma_c(M;\II^{\infty,\bullet}_{M,M \setminus U}),
$$
where the first isomorphism follows from the definition of proper direct image and the second one follows applying $R\Gamma_c(M;\cdot)$ to Theorem \ref{WDR}.\qed

\begin{oss} As a consequence of Corollary \ref{SDR}, we obtain a result of \cite{AG10}: the cohomology of the de Rham complex with Schwartz coefficients $H^\bullet DR_\S(M)$ of a Nash manifold $M$ is isomorphic to the compact support cohomology $H_c^\bullet(M)$ of $M$.  We first prove it for $\R^n$. We can see $\R^n$ as an open subset of $\mathbb{S}^n$. Let $U$ be an open semialgebraic subset of  $\R^n$. Then
$$
H^\bullet_c(U) \simeq H^\bullet(\mathbb{S}^n;\II^{\infty,\bullet}_{\mathbb{S}^n,\mathbb{S}^n \setminus U}) \simeq H^\bullet DR_\S(U).
$$
Using the fact that a Nash manifold has a finite cover consisting of open submanifolds Nash diffeomorphic to $\R^n$ one obtains the result.
\end{oss}

\section{A CDGA quasi-isomorphism}

Let $Z$ be a closed subanalytic subset of $M$ and consider the de Rham complex of sheaves
$$
\CC_Z \wtens \C^{\infty,\bullet}_M.
$$
%It is defined by the exact sequence of (de Rham) complexes
%$$
%\exs{\CC_U \wtens \C^{\infty,\bullet}_M}{\C^{\infty,\bullet}_M}{\CC_Z \wtens \C^{\infty,\bullet}_M}
%$$
%Note that this definition goes back to a paper of Kashiwara-Schapira of '96 and the first one is the complex with coefficients vanishing with all their derivatives outside $U$.
%$\cdot \wtens \C^\infty_M$ is an exact covariant functor defined in the category of $\R$-constructible sheaves. \\
Let us consider a resolution of $Z$
$$
0 \to \oplus \CC_{U_1} \to \cdots \to \oplus \CC_{U_k} \to \CC_Z \to 0
$$
made by subanalytic open subsets (we may also assume they are contractible and locally cohomologically trivial), then $\CC_Z \wtens \C^{\infty,\bullet}_M$ is defined by a sequence of chain complexes whose entries are of this kind
$$
\CC_V \wtens \C^{\infty,\bullet}_M,
$$
with $V$ locally cohomologically trivial. Now, by definition of $\CWM$, we have $\CC_V \wtens \C^{\infty}_M \simeq \imin\rho\ho(\CC_{\overline{V}},\CWM)$ when $V$ is locally cohomologically trivial. It follows from Lemma \ref{deRham} that the chain complex of (subanalytic) sheaves
$$
\CC_M \to \C^{\infty,{\rm w},\bullet}_M
$$
is exact. Moreover $\CC_M$ and $\C^{\infty,{\rm w}}_M$ are acyclic with respect to $\ho(\CC_{\overline{V}},\cdot)$ ($V$ locally cohomologically trivial) and $\imin\rho$ is exact, hence the chain complex
$$
\CC_V \simeq \ho(\CC_{\overline{V}},\CC_M) \to \imin\rho\ho(\CC_{\overline{V}},\C^{\infty,{\rm w},\bullet}_M) \simeq \CC_V \wtens \C^{\infty,\bullet}_M
$$
is exact. This implies that the chain complex
$$
\CC_Z \to \CC_Z \wtens \C^{\infty,\bullet}_M
$$
is exact as well.
%Moreover $F\wtens\C^{\infty,\bullet}_M$ is a complex of sheaves which are injective with respect to the functor $\Gamma(W,\cdot)$, $W$ open in $M$.
Let us consider the exact sequence of chain complexes
$$
0 \to \CC_U \wtens \C^{\infty,\bullet}_M \to \C^{\infty,\bullet}_M \to \CC_Z \wtens \C^{\infty,\bullet}_M \to 0
$$
and apply the functor $\Gamma(W;\cdot)$. Being the Whitney tensor product a soft sheaf,
%(for any $\R$-constructible sheaf),
we get an exact sequence of chain complexes
$$
0 \to \Gamma(W;\CC_U \wtens \C^{\infty,\bullet}_M) \to \Gamma(W;\C^{\infty,\bullet}_M) \overset{J}{\to} \Gamma(W;\CC_Z \wtens \C^{\infty,\bullet}_M) \to 0.
$$
Remark that $J$ is a CDGA (commutative differential graded algebra) morphism.
We have
\begin{eqnarray*}
\Gamma(W;\CC_U \wtens \C^{\infty,\bullet}_M) & \simeq & R\Gamma(W;\CC_U), \\
\Gamma(W;\C^{\infty,\bullet}_M) & \simeq & R\Gamma(W;\CC_M), \\
\Gamma(W;\CC_Z \wtens \C^{\infty,\bullet}_M) & \simeq & R\Gamma(W;\CC_Z).
\end{eqnarray*}
Suppose that $Z$ is a deformation retract of $W$. Then
\begin{eqnarray*}
R\Gamma(W;\CC_U) & \simeq & 0, \\
R\Gamma(W;\CC_M) & \simeq & R\Gamma(Z,\CC_Z), \\
R\Gamma(W;\CC_Z) & \simeq & R\Gamma(Z;\CC_Z),
\end{eqnarray*}
%The left-hand side represents the complex $R\Gamma(W;\CC_M) \simeq R\Gamma(Z;\CC_Z)$,
%being $Z$ a deformation retract of $W$.
%The right-hand side represents the complex $R\Gamma(M;\CC_Z) \simeq R\Gamma(Z;\CC_Z)$.
Hence for such a $W$
$$
\Gamma(W;\C^{\infty,\bullet}_M) \overset{J}{\to} \Gamma(W;\CC_Z \wtens \C^{\infty,\bullet}_M).
$$
is a quasi-isomorphism of CDGA.

\begin{oss} \label{HZHWZ} Since every neighborhood of a subanalytic closed subset $Z$ of $M$ contains a neighborhood $W$ which has a deformation retract to $Z$, we get a quasi-isomorphism of CDGA
$$
\Gamma(M;\CC_Z \otimes \C^{\infty,\bullet}_M) \simeq \lind W \Gamma(W;\C^{\infty,\bullet}_M) \overset{J}{\to} \lind W \Gamma(W;\CC_Z \wtens \C^{\infty,\bullet}_M) \simeq \Gamma(M;\CC_Z \wtens \C^{\infty,\bullet}_M)
$$
where $W$ ranges through the family of neighborhoods having a deformation retract to $Z$. Thanks to this fact one can easily prove the quasi-isomorphism of sheaves of CDGA
$$
\CC_Z \otimes \C^{\infty,\bullet}_M \iso \CC_Z \wtens \C^{\infty,\bullet}_M.
$$
\end{oss}

\begin{oss} Thanks to the exactness of $\cdot\otimes\C^\infty_M$ and $\cdot\wtens\C^\infty_M$, one can extend the quasi-isomorphism of Remark \ref{HZHWZ} to the case of an open subanalytic subset $U$ and then, inductively, to the case of a $\R$-constructible sheaf $F$. Namely, we have a quasi-isomorphism of sheaves of CDGA
$$
F \otimes \C^{\infty,\bullet}_M \iso F \wtens \C^{\infty,\bullet}_M.
$$
\end{oss}

\begin{oss} As a further application, one can extend the results of \cite{BP15} (as pointed out by the authors in the introduction) to the case of subanalytic sets. One first remarks (following their notations) that the Whitney-de Rham cohomology of $Z$ is isomorphic as a CGA to its singular cohomology with real coefficients:
$$
H^\bullet(Z,\R) \simeq H^\bullet_W(Z),
$$
and these isomorphisms are compatible with the CGA structures. Then, when the set $Z$ is simply connected, one can prove that the de Rham complex with Whitney coefficients $\Omega^\bullet_W(Z)$ determines the real homotopy type of $Z$. This also implies that the Hochschild homology of the differential graded algebra $\Omega^\bullet_W(Z)$ is isomorphic to the cohomology of the free loop space $\mathcal{L}Z$.
\end{oss}

\end{document}